\documentclass[a4paper,12pt]{article}

\usepackage{amsmath,amssymb,amsfonts,epsfig,cite,setspace}
\usepackage[paper=letterpaper,margin=1.2in]{geometry}

\begin{document}


\renewcommand{\thefootnote}{\fnsymbol{footnote}}
\centerline{\Huge Polynomial Roots and Calabi-Yau Geometries}
~\\
\vskip 2mm
\centerline{
{\large Yang-Hui He}\footnote{\tt hey@maths.ox.ac.uk}
}
~\\
{\scriptsize
\begin{center}
\begin{tabular}{ll}
  $^1$ & {\it School of Physics, NanKai University,}
  {\it Tianjin, 300071, P.R.~China}\\
  $^2$ & {\it Department of Mathematics, City University London,} 
  {\it Northampton Square, London EC1V 0HB, U.K.}\\
  $^3$ & {\it Merton College, University of Oxford,}
  {\it OX1 4JD, U.K.}\\
\end{tabular}
\end{center}
}

~\\
~\\

\newcommand{\todo}[1]{{\bf ?????!!!! #1 ?????!!!!}\marginpar{$\Longleftarrow$}}
\newcommand{\eref}[1]{Eq.~(\ref{#1})}
\newcommand{\sref}[1]{\S\ref{#1}}
\newcommand{\fref}[1]{Figure~\ref{#1}}
\newcommand{\nn}{\nonumber}
\newcommand{\comment}[1]{}

\newcommand{\CM}{{\cal M}}
\newcommand{\CN}{{\cal N}}
\newcommand{\CZ}{{\cal Z}}
\newcommand{\cO}{{\cal O}}
\newcommand{\cB}{{\cal B}}
\newcommand{\cC}{{\cal C}}
\newcommand{\cD}{{\cal D}}
\newcommand{\cE}{{\cal E}}
\newcommand{\cF}{{\cal F}}
\newcommand{\cX}{{\cal X}}
\newcommand{\IA}{\mathbb{A}}
\newcommand{\IP}{\mathbb{P}}
\newcommand{\IR}{\mathbb{R}}
\newcommand{\IC}{\mathbb{C}}
\newcommand{\IF}{\mathbb{F}}
\newcommand{\IV}{\mathbb{V}}
\newcommand{\II}{\mathbb{I}}
\newcommand{\IZ}{\mathbb{Z}}
\newcommand{\re}{{\rm Re}}
\newcommand{\im}{{\rm Im}}
\newcommand{\li}{{\rm Li}}

\newcommand{\tmat}[1]{{\tiny \left(\begin{matrix} #1 \end{matrix}\right)}}

\newcommand{\diff}[2]{\frac{\partial #1}{\partial #2}}

\newtheorem{theorem}{\bf THEOREM}
\def\thetheorem{\thesection.\arabic{theorem}}
\newtheorem{conjecture}{\bf CONJECTURE}
\def\thetheorem{\thesection.\arabic{conjecture}}
\newtheorem{observation}{\bf OBSERVATION}
\def\thetheorem{\thesection.\arabic{observation}}

\def\theequation{\thesection.\arabic{equation}}
\newcommand{\setall}{\setcounter{equation}{0}
        \setcounter{theorem}{0}}
\newcommand{\setequation}{\setcounter{equation}{0}}

\begin{abstract}
The examination of roots of constrained polynomials dates back at least to Waring and to Littlewood.
However, such delicate structures as fractals and holes have only recently been found.
We study the space of roots to certain integer polynomials arising naturally in the context of Calabi-Yau spaces, notably Poincar\'e and Newton polynomials, and observe various salient features and geometrical patterns.
\end{abstract}

\newpage
\tableofcontents

\vspace{2cm}

\section{Introduction and Summary}\setall
The subject of roots of mono-variate polynomials is, without doubt, 
an antiquate one, and
has germinated an abundance of fruitful research over the ages.
It is, therefore, perhaps surprising that any new statements could at all be
made regarding such roots.
The advent of computer algebra, chaotic phenomena, and  random ensembles
has, however, indeed shed new light upon so ancient a metier.

Polynomials with constrained coefficients and form, though permitted to
vary randomly, have constituted a vast field itself.
As far back as 1782, Edward Waring, in relation to his famous problem
on power summands, had shown that for cubic polynomials with random real
coefficients, the ratio of the probability of finding non-real zeros versus
that of not finding non-real zeros is less than or equal to 2.
Constraining the coefficients to be integers within a fixed range has, too,
its own history.
It was realised in \cite{BlockPolya} that a degree $n$ random polynomial 
$P(z) = 1 + \sum\limits_{k=1}^n a_k z^k$ with $a_k = -1,0,1$ distributed
evenly, the expected number $\nu_n$ of real roots is of order $\cO(n^{1/2})$
asymptotically in $n$.
This was furthered by \cite{LittlewoodOfford} to be essentially independent of
the statistics, in that $\nu_n$ has the same asymptotics\footnote{Cf.~also \cite{Hannay,PeresVirag}.}, as much for 
$a_k$ being evenly distributed real numbers, in $[-1,1]$, or as Gaussian distributed in $(-\infty,\infty)$.

Continual development ensued\footnote{q.v.~also \cite{randompoly}.}, notably by Littlewood \cite{Littlewood}, Erd\H{o}s \cite{Erdos}, Hammersley \cite{Hammersley} and Kac \cite{Kac}.
Indeed, a polynomial with coefficients only taking values as $\pm 1$ has come to be known as a {\bf Littlewood polynomial} and the {\it Littlewood Problem} asks for the the precise asymptotics, in the degree, of such polynomials taking values, with complex arguments, on the unit circle.
The classic work of Montgomery \cite{Montgomery} and Odlyzko \cite{Odlyzko},
constituting one of the most famous computer experiments in mathematics\footnote{q.~v.~Section 3.1 of \cite{He:2010jh} for some recent remarks on the distributions.},
empirically showed that the distribution of the (normalized) spacings between
successive critical zeros of the Riemann zeta function is the same as that of
a Gaussian unitary ensemble of random matrices, whereby infusing our subject with issues of uttermost importance.

\begin{figure}[!h!t!b]
\centerline{
(a) \includegraphics[trim=0mm 0mm 0mm 0mm, clip, width=2.0in]{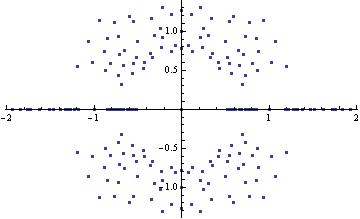}
(b) \includegraphics[trim=0mm 0mm 0mm 0mm, clip, width=2.0in]{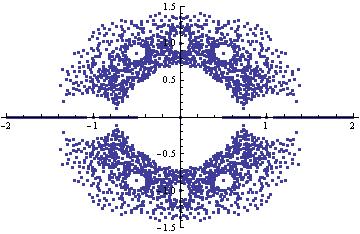}
(c) \includegraphics[trim=0mm 0mm 0mm 0mm, clip, width=2.0in]{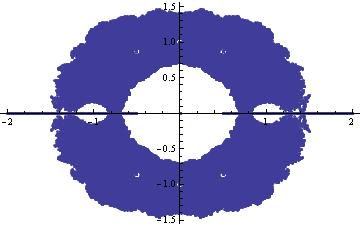}
}
\caption{{\sf 
The position, on the complex plane, of the zeros of 50000 random integer polynomials with coeffcients
$-1$, 0 or 1, for degrees upto 4, in part (a), 6, in part (b) and 10, in part (c).
}
\label{f:randRange-1to1}}
\end{figure}

Subsequently, combining the investigation of zeros and of random polynomials, Odlyzko and Poonen \cite{OdPoo} studied the zeros of Littlewood-type polynomials by setting the coefficients to 0 and 1; they provided certain bounds as well as found interesting fractal structures.
Thus inspired and with the rapid advance of computational power, Borwein et al.~constructed various plots of zeros of constrained random polynomials and many remarkable features were instantly visible \cite{BorweinJorg,BBCGLM}.
We can readily demonstrate this with Mathematica$^{{\tiny \textregistered}}$ \cite{math}, as is shown in \fref{f:randRange-1to1}.
In the figure, we take a sample of 50000 random polynomials with coefficients $-1$, 0 or 1 up to various degrees, and plot, on the complex plane, their zeros.
Not only do we see fractal behaviour\footnote[1]{Cf.~discussions in \cite{Franco:2004jz,Hannay} on chaotic dualities in field theories} near the boundaries, the nature of the holes are intimately related \cite{BBBP} to the Lehmer-Mahler Conjecture: that the Mahler measure
$M(P) := \exp\left( \frac{1}{2\pi} \int_0^{2\pi} \log \left| P(e^{i \theta}) \right| d\theta \right)$ of any integral polynomial $P(z)$ (which is not a multiple of cyclotomic polynomials) should be bounded below by that of $z^{10}-z^9+z^7-z^6+z^5-z^4+z^3-z+1$, which is approximately $1.17$.

High resolution variants of \fref{f:randRange-1to1} have been considered recently by Christensen \cite{Christensen}, J\"orgenson \cite{Jorgenson} and Derbyshire \cite{Derbyshire}, inter alia, and many beautiful pictures can be found\footnote{Cf.~also a nice account in \cite{baez}.}. Particular striking are the coloured density plots in \cite{Derbyshire}.

An interesting query, in somewhat reverse direction to the above line of thought, was posed in \cite{Lin}: recalling that the Lee-Yang Circle Theorem placed severe constraints on the generating function of the partition function of the Ising Model, the author asked if one could statistically test whether a given Laurent polynomial could, in fact, be the Jones polynomial of a knot.
Using a landscape of knots generated by the programme ``knotscape'' \cite{knotscape}, the said work investigated many distribution properties of the zeros of known Jones polynomials.

\vspace{1cm}

Enchanted by this motif which has threaded varying developments over the decades while persistently generating new perspectives, a question immediately springs to mind.
One of the central topics of both modern mathematics and theoretical physics is undoubtedly that of Calabi-Yau geometries.
A key feature is the super-abundance thereof.
In complex dimension one, there is only the torus; in dimension two, there are the 4-torus and the K3-surface; however, for dimension three and above, no classification is known and already a plethora has been constructed.
The first database was that of so called CICYs, or, complete intersection Calabi-Yau threefolds in products of projective spaces \cite{Candelas:1987kf} as well as hypersurfaces in weighted $\IC\IP^4$ \cite{Candelas:1989hd}, by Candelas et al. Then, over more than a decade, Kreuzer and Skarke formulated and compiled an impressive list of on the order of $10^{10}$ threefolds as hypersurfaces in toric varieties \cite{Kreuzer:2000xy}.
Finding new patterns in this vast distribution of manifolds has seen some recent activity \cite{Candelas:2007ac,Candelas:2008wb,Braun:2010vc}.
Indeed, the multitude of these geometries is at the core of the so-called vacuum degeneracy problem of superstring theory and constitutes a part of the landscape issue \cite{Denef:2006ad}.

Along a parallel vein, the space of non-compact (singular) Calabi-Yau spaces (as affine varieties) has also been extensively explored, notably by Hanany et al.~over the past decade \cite{Hanany:1998sd,Feng:2000mi,Hanany:2005ve,Franco:2005sm,Forcella:2008bb,Hanany:2008gx,Davey:2009et,Davey:2010px}, especially those which admit a toric description \cite{Feng:2000mi}; the discovery of their intimate relation to dimer models and brane tilings \cite{Hanany:2005ve,Franco:2005sm,Davey:2009et} has also generated some excitement.
These supplant yet another corner in the landscape of geometries and associated supersymmetric vacua.

\vspace{1cm}

Thus motivated, many tasks lend themselves to automatic investigation; we here give a precis of some key points.
In Section 2, we begin with the compact, smooth Calabi-Yau manifolds.
As mentioned above, there had been much effort in classifying and constructing these, especially in complex dimensions three and four.
An immediate polynomial, of constrained form and integer coefficients, and yet succinctly encoding some topological information, is the Poincar\'e polynomial, which can be readily written in terms of the Hodge numbers.
We take the ``experimental data'' of all the known Hodge numbers of the threefolds and fourfolds spanning two decades of work and plot the complex roots of the associated Poincar\'e polynomials in \fref{f:h11h12CY3} and \fref{f:poincareCY4} respectively.

Much intricate structures are clearly visible.
These are then contrasted with a ``standard background sample'', namely the roots of random integer sextic and octic polynomials, with unit leading coefficient and vanishing linear term, drawn in \fref{f:palin} and \fref{f:palin4}.
From such collections are extracted the sub-class of those which admit $-1$ as a roots, which, by a theorem from differential geometry, correspond to spaces which have more than one isometry.
Interestingly, they correpond to self-mirror threefolds and ``quasi''-self-mirror fourfolds.
We plot the roots for these in \fref{f:selfmirror} and \fref{f:selfmirrorCY4},and see that they furnish certain substrata of the conglomerate plots mentioned above.

Thenceforth we move on to non-compact Calabi-Yau geometries in Section 3.
There, too, is a plenitude of examples, most notably those which are toric.
We focus on toric threefolds because these have planar toric diagrams as lattice points in $\IZ^2$ due to the Calabi-Yau condition.
Once more, a natural polynomial invites itself: the Newton polynomial.
The Riemann surface corresponding to this bi-variate polynomial has been a central subject to the gauge and brane theories in the context of string theory.
Moreover, the two important projections thereof, viz., the amoeba and alga projections, have provided many beautiful Monte Carlo plots, illustrating deep algebraic geometry as well as gauge theory.
Because we are confronted with two complex variables, we need to slightly deviate from our theme of complex roots; instead, we find it expedience to regard the variables as real and consider the real projection of the Riemann surface.
Subsequently, we can study the ensemble of real turning (critical) points of these planar curves.

Again, we resort to ``actual data'' and focus on the most well-known affine toric threefold geometries - as shown in \fref{f:toricdiag} - corresponding to local Calabi-Yau singularities, including, of course, the famous conifold.
To each space, we find the collective of critical points in $\IR^2$ as we vary the integer coefficients - commonly known as multiplicities in the dimer model interpretation - of the Newton polynomial, and plot them in \fref{f:critical}.
We see a sensitive dependence of the emergent subtle structures upon the choice of toric data.

In many respects we have taken a very pragmatic and empirical approach toward the data accumulated over many years of theoretical research, of quantity large enough to justify experimentation.
To this philosophy of ``experimental mathematics'' we adhere throughout, observe wherever we should and infer wherever we may.
Without much to do, therefore, let us delve into the details of the issues summarized above.

\vspace{0.7cm}

{\small
\begin{singlespace}
{\it To the occasion of the happy Installation of Professor Sir Martin Taylor as the Fiftieth Warden of Merton College, Oxford, on the Second Day of the Month of October, in the Year of Our Lord Two Thousand and Ten, and to the noble retirement of Professor Dame Jessica Rawson, this humble brief note is dedicated. 

Vivat Custos, vivat Collegium, \& Stet Fortuna Domus!
}
\end{singlespace}
}

\section{Compact Calabi-Yau Manifolds, Poincar\'e Polynomials and Complex Roots}\setall
An important quantifying polynomial for a smooth compact manifold $X$ is the Poincar\'e polynomial, which is a generating function for topological invariants of $X$ (say of dimension $n$) :
\begin{equation}\label{PoincareP}
P(t;X) = \sum\limits_{i=0}^n b_i t^i \ ,
\end{equation}
with the $b_i$ being the $i$-th Betti number.
Indeed, this seems a more natural candidate for our present studies than some because other famous polynomials such as the Hilbert polynomial or the numerator of the Hilbert series (which of late have been instrumental in counting BPS operators \cite{Benvenuti:2006qr,Forcella:2009vw,Hanany:2008sb}) are not topological invariants and depend on the specific projective embedding.
Furthermore, by definition, at $t=-1$, the polynomial evaluates to the Euler characteristic; this will be of significance shortly.

The zeros of the Poincar\'e polynomial have rather remarkable properties.
It was conjectured that \cite{Emery} that if the rank of the manifold $X$ is greater than 1, where rank is defined to the the maximal number of everywhere independent, mutually commuting, vector fields on $X$, i.e., the number of isometries, then $-1$ is a multiple root of the Poincar\'e polynomial of $X$.
Unfortunately this conjecture was shown to be false \cite{Brendon}.
Nevertheless, it still holds that the rank of $X$ exceeds unity if and only if $-1$ is a multiple root of $P(t;X)$.

Moreover, of number theoretic and arithmo-geometric significance is the fact that certain Poincar\'e polynomials exhibit Riemann Hypothesis behaviour \cite{P-arith}, in analogy to the the Hasse-Weil zeta local zeta functions.
Recently, alignment of zeros of Hilbert polynomials have been studied by \cite{Rodriguez} in relation to zeta functions.

\subsection{Calabi-Yau Threefolds}
Our focus will be on \eref{PoincareP}.
First, let us study the case of Calabi-Yau threefolds, which have been of the greatest interest, at least historically. 
Because we are dealing with complex (K\"ahler) manifolds, Hodge decomposition implies that $b_i(X) = \sum\limits_{p,q} h^{p,q}(X)$, with $h^{p,q}(X) = \dim H_{\bar{\partial}}^{p,q}(X)$ the dimensions of the Dolbeault cohomology groups. Indeed, for (compact, smooth, connected) Calabi-Yau threefolds, the Hodge diamond, and subsequently the Betti numbers and the Poincar\'e polynomial, can be written as:
\begin{equation}
\begin{array}{lllr}
\begin{array}{ccccccc}
& & & 1 & & &  \\
& & 0 &  & 0 & & \\
& 0 & &  h^{1,1} & & 0 &  \\
1 & & h^{2,1} & & h^{2,1} & &  1 \\
& 0 & &  h^{1,1} & & 0 &   \\
& & 0 &  & 0 & &   \\
& & & 1 & & &   \\
\end{array}
&&&
\begin{array}{l}
(b_0=b_6,b_1=b_5,b_2=b_4,b_3) = (1,0,h^{1,1},2+2h^{2,1}); \\ ~\\ ~\\
P(t;X) = 1 + h^{1,1} t^2 + (2+2h^{2,1}) t^3 + h^{1,1} t^4 + t^6 \ .
\end{array}
\end{array}
\label{PCY3}
\end{equation}
That the Poincar\'e polynomial is palindromic is obvious and follows from Poincar\'e duality.

\begin{figure}[!h!t!b]
\centerline{
(a) \includegraphics[trim=0mm 0mm 0mm 0mm, clip, width=2.0in]{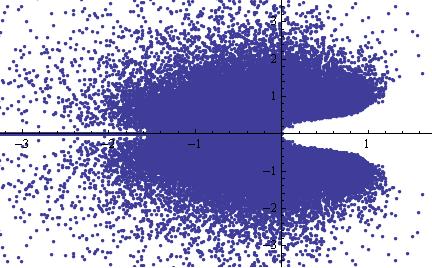}
(b) \includegraphics[trim=0mm 0mm 0mm 0mm, clip, width=2.0in]{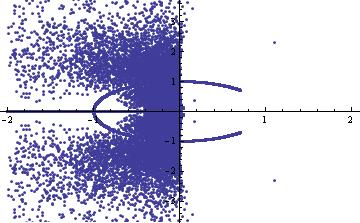}
(c) \includegraphics[trim=0mm 0mm 0mm 0mm, clip, width=2.0in]{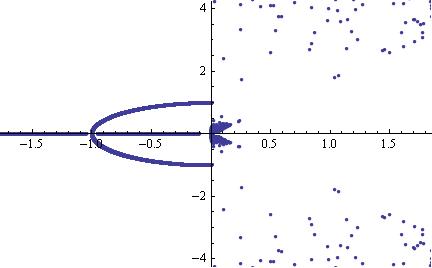}
}
\caption{{\sf 
(a) The position, on the complex plane, of the zeros of 50000 random integer degree six polynomials with coefficients between 0 and 1000.
(b) The same, but with monic palindromic sextics.
(c) Monic palindromic sextics, and with linear term vanishing.
}
\label{f:palin}}
\end{figure}

Therefore, our first constraint is palindromicity to which we shall presently restrict. We recall the roots of a completely random sample of integer polynomials with coefficients in $[-1,1]$ up to the sextic in part (b) of \fref{f:randRange-1to1}.
In \fref{f:palin}, we plot, in part (a), a sample of 50000 random integer sextic polynomials with coefficients in $[0,1000]$ (making sure that the highest coefficient at degree 6 is not 0) as a comparative norm.
Next, in part (b), we plot the same, but for monic palindromic sextics, i.e., 
$P(t) = 1 + b_1 t + b_2 t^2 + b_3 t^3 + b_2 t^4 + b_1 t^5 + t^6$. Then, in (c), we restrict once more, with some foresight, so that the linear term vanishes, i.e., $P(t) = 1 + b_2 t^2 + b_3 t^3 + b_2 t^4 + t^6$.
We see that upon the condition of palindromicity, there is a marked emergence of roots on the unit circle; this of course arises from the symmetric terms of the form $e^{it}$ combining to give (co)sines whose reality then facilitates the addition to zero.
The symmetry about the $x$ axis is simply that all roots appear in conjugate pairs because our polynomials have real coefficients.

\begin{figure}[!h!t!b]
\centerline{
(a) \includegraphics[trim=0mm 0mm 0mm 0mm, clip, width=3.0in]{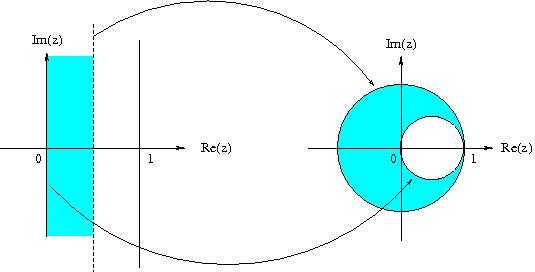}
(b) \includegraphics[trim=0mm 0mm 0mm 0mm, clip, width=3.0in]{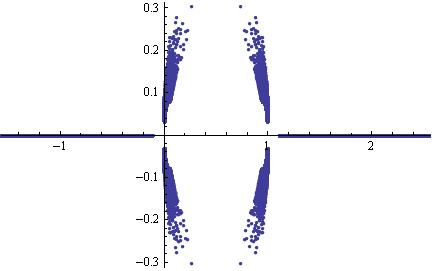}
}
\caption{{\sf 
(a) The conformal map $z \to \frac{z}{z-1}$ takes the left half of the critical strip to the inside of the unit disk, with the boundaries mapped as shown by the arrows. It takes the mirror image, in the right half of the critical strip, to the complement of the unit disk. The inverse map is given by $z \to \frac{z}{z+1}$.
(b) The position of 50000 randomly integer monic palindromic sextic polynomials with vanishing linear/quintic terms and with coefficients ranging in $[0,1000]$, applying apply the map $z \to \frac{z}{z+1}$.
}
\label{f:conformal}}
\end{figure}

For our amusement, seeing the form of the semi-unit-circular shape being prominent, we are reminded of the conformal map $z \to \frac{z}{z-1}$ which takes the unit circle to the critical strip of the Riemann Hypothesis, as shown in detail by part (a) of \fref{f:conformal}.
We take the space of roots of monic palindromic sextics with vanishing linear/quintic terms from part (c) of \fref{f:palin}, apply the inverse map $z \to \frac{z}{z+1}$ to map to the critical strip and re-do the plot in part (b) of \fref{f:conformal}.
We see that the plot resembles the zero-free region of the Riemann zeta function inside the critical strip.

Founded upon these above discussions, we are ready to approach ``actual data''.
We now collect all known Calabi-Yau threefolds, which come from three major databases, viz., the aforementioned CICYs, the hypersurfaces in toric varieties, as well as the collection of individually tailored ones of small Hodge numbers (cf.~\cite{Candelas:2007ac,Candelas:2008wb}).
These total, respectively, 30108, 266 and 54 distinct pairs of Hodge numbers $(h^{1,1}(X),h^{2,1}(X))$.
In all, there are 30237 distinct pairs (of course, each with much degeneracy) of Hodge numbers; to our present knowledge, these are all the ones circulated in the literature.

\begin{figure}[!h!t!b]
\centerline{
(a) \includegraphics[trim=0mm 0mm 0mm 0mm, clip, width=2.0in]{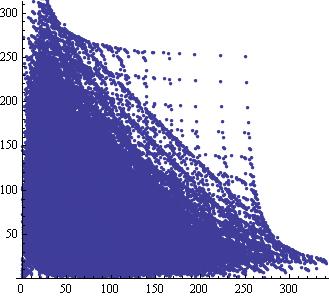}
(b) \includegraphics[trim=0mm 0mm 0mm 0mm, clip, width=4.0in]{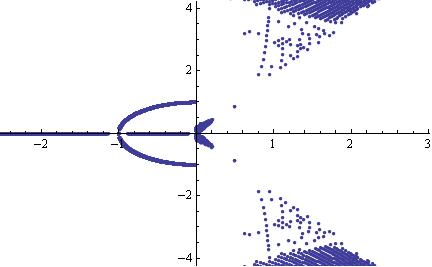}
}
\caption{{\sf 
(a) The Hodge numbers, with $h^{1,1}$ as the abscissa and $h^{2,1}$, the ordinate, of all the known Calabi-Yau threefolds.
(b) The position, on the complex plane, of the zeros of their Poincar\'e polynomials.
}
\label{f:h11h12CY3}}
\end{figure}

We plot\footnote{Traditionally, the now-famous plot, first appearing in \cite{Candelas:1989hd}, is done with $\chi = 2(h^{1,1}-h^{2,1})$ as abscissa and $h^{1,1}+h^{2,1}$ as ordinate.} these, with $h^{1,1}$ as the abscissa and $h^{2,1}$, ordinate, in part (a) of \fref{f:h11h12CY3}, the largest amongst these is $(491,11)$. 
Note that because of mirror symmetry, there is a symmetry interchanging the two co\"ordinates.
It is still an open question whether there exists any Calab-Yau threefold whose either Hodge number exceeds 491, a bound which has defied constructions so far.
This is why in our random standard background sample in \fref{f:palin}, we have conveniently selected the largest integer coefficient to be $1000 \sim 2 \cdot (491+1)$. 
In part (b) of the said Figure, we plot, on the complex plane, the roots of the Poincar\'e polynomials of all these known threefolds.
Because we are dealing with polynomials of non-negative coefficients, there should be many generic roots with negative real parts.
Comparing with the random sample in part (c) of \fref{f:palin}, we see a beautiful clustering of points in the first quadrant (and, by complex conjugation, the fourth).

\begin{figure}[!h!t!b]
\centerline{
(a) \includegraphics[trim=0mm 0mm 0mm 0mm, clip, width=2.0in]{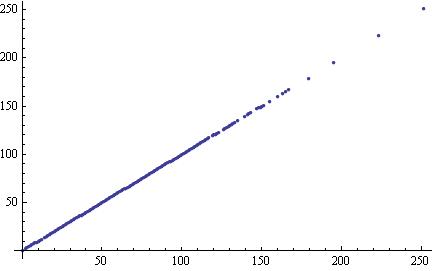}
(b) \includegraphics[trim=0mm 0mm 0mm 0mm, clip, width=3.0in]{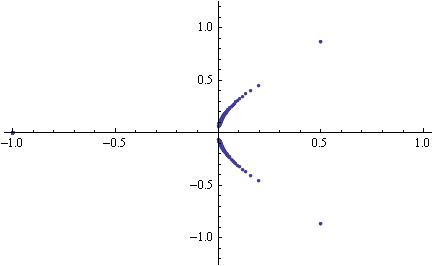}
}
\caption{{\sf 
(a) The Hodge number of self-mirror Calabi-Yau threefolds; these have $-1$ as a root of the Poincar\'e polynomial.
(b) All of the roots of the Poincar\'e polynomial of these self-mirrors.
}
\label{f:selfmirror}}
\end{figure}

Next, let us test for how many Poincar\'e polynomials $-1$ is a root; these, as mentioned above, would correspond to manifolds which have more than one isometries.
Interestingly, of the some 30000, only 148 pass the test.
These turn out to be only the 148 known self-mirror manifolds; we plot their Hodge numbers in Part (a) of \fref{f:selfmirror} (the values of Hodge numbers range from $(1,1)$ to $(251,251)$, skipping many high values, as well as the number 13).
Indeed, this is a simple consequence of palindromicity as one sees that, upon substituting $t=-1$ into $P(t;X)$ in \eref{PCY3}, we obtain $P(-1;X) = 2(h^{1,1}-h^{2,1}) = \chi(X)$, the Euler characteristic.
In Part (b), we plot all the other roots as well, and see that these constitute a portion of the small crescent around the origin.

\subsection{Calabi-Yau Fourfolds}
Having explored the space of Calabi-Yau threefolds, it is automatic to proceed onto the space of fourfolds, a relative {\it terra incognita}.
We again resort to the wonderful database compiled by Kreuzer-Skarke \cite{Kreuzer:2000xy,Kreuzer:1997zg}.
Now, there are, totaling the hypersurfaces in toric fivefolds, 14598161 manifolds, with 3015056 distinct triplet of Hodge numbers.
To explain this triplet notation, we remind the reader of the Hodge diamond of compact, connected, smooth Calabi-Yau fourfolds, adhering to the nomenclature and explanation of \cite{Lynker:1998pb}:
\begin{equation}
\begin{array}{cccc}
\begin{array}{ccccccccc}
 &&&& 1 &&&& \\
  &&&0 && 0 &&& \\ 
  &   & 0 & & h^{1,1} & & 0 & & \\
  & 0 & & h^{2,1} & &  h^{2,1}  & & 0 & \\
1 &   & h^{3,1} & & h^{2,2} & &  h^{3,1} & & 1\\
  & 0 & & h^{2,1} & &  h^{2,1}  & & 0 & \\
  &   & 0 & & h^{1,1} & & 0 & & \\
  &&&0 && 0 &&& \\
 &&&& 1 &&&& \\ 
\end{array}
&&&
\begin{array}{l}
h^{2,2} = 44 + 4 h^{1,1} - 2 h^{2,1} + 4 h^{3,1};\\
~\\~\\
(b_0=b_8,\; b_1=b_7,\; b_2=b_6,\;b_3=b_5,\;b_4) = \\
(1,\; 0, ; h^{1,1}, \; 2h^{2,1}, \; 2 + 2h^{3,1} + h^{2,2}) = \\
(1,\;0, \;h^{1,1},\; 2h^{2,1}, \; 46+4h^{1,1}-2h^{2,1}+6h^{3,1})
\end{array}
\end{array}
\label{bCY4}
\end{equation}
We note that though seemingly there are four degrees of freedom, owing to topological constraints in complex dimension four or higher, as exhibited by the above relation of $h^{2,2}$ with the others, there are really only three independent Hodge numbers, which we choose as $(h^{1,1}, h^{2,1}, h^{3,1})$; in terms of this triplet we express the Betti numbers, as shown above.
Consequently, we can write the Poincar\'e polynomial of the fourfold, from the Betti numbers in \eref{bCY4}, as
\begin{equation}
P(t;X) = 1 + h^{1,1} t^2 + 2h^{2,1} t^3 + (46+4h^{1,1}-2h^{2,1}+6h^{3,1}) t^4
+ 2h^{2,1} t^5 + h^{1,1} t^6 + t^8 \ .
\label{PCY4}
\end{equation}

\begin{figure}[!h!t!b]
\centerline{
(a) \includegraphics[trim=0mm 0mm 0mm 0mm, clip, width=3.5in]{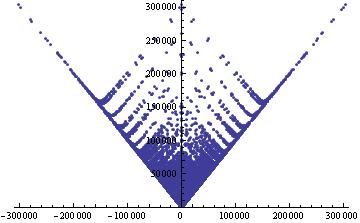}
(b) \includegraphics[trim=0mm 0mm 0mm 0mm, clip, width=2.5in]{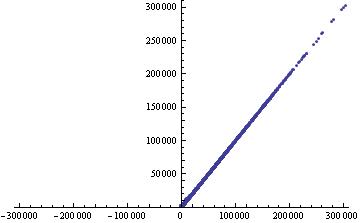}
}
\caption{{\sf 
(a) The Hodge numbers, with $h^{1,1}-h^{3,1}$ as the abscissa and $h^{1,1}+h^{3,1}$, the ordinate, of the fourfolds from Kreuzer-Skarke's database.
(b) The same, but with $h^{1,1}-h^{2,1}$ as the abscissa and $h^{1,1}+h^{2,1}$ as the ordinate.
}
\label{f:CY4-plot}}
\end{figure}

Following \cite{Lynker:1998pb}, we plot $h^{1,1}+h^{3,1}$ as ordinate versus $h^{1,1}-h^{3,1}$ as abscissa, which demonstrates mirror-like behaviour\footnote{Though in {\it cit.~ibid.}, only the hypersurfaces in weighted $\IC\IP^5$ were considered, whereas here we plot the entirety of the known fourfolds.}.
We also plot $h^{1,1}+h^{2,1}$ versus $h^{1,1}-h^{2,1}$, showing that the behaviour in the applicate direction is rather trivial.
These are shown in \fref{f:CY4-plot}.

\begin{figure}[!h!t!b]
\centerline{
(a) \includegraphics[trim=0mm 0mm 0mm 0mm, clip, width=2.0in]{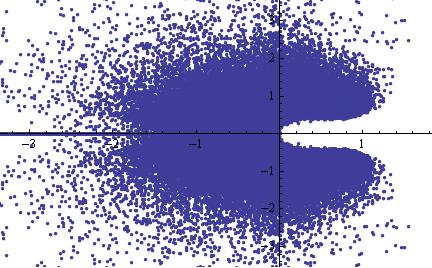}
(b) \includegraphics[trim=0mm 0mm 0mm 0mm, clip, width=2.0in]{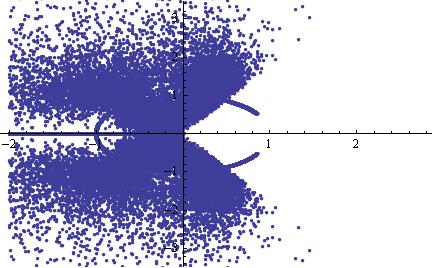}
(c) \includegraphics[trim=0mm 0mm 0mm 0mm, clip, width=2.0in]{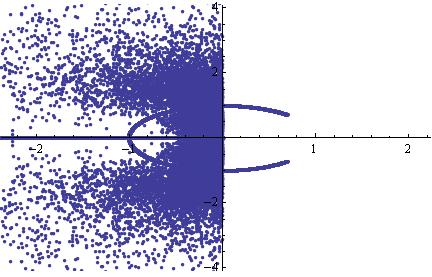}
}
\caption{{\sf 
(a) The position, on the complex plane, of the zeros of 50000 random integer polynomials with coeffcients between 0 and 2500000.
(b) The same, but with monic palindromic octics.
(c) Monic palindromic octics, and with linear term vanishing.
}
\label{f:palin4}}
\end{figure}

We now repeat the experiment undertaken for threefolds.
First, we plot the space of generic roots, and present them in \fref{f:palin4}.
In part (a), a sample of 50000 random integer octic polynomials with coefficients in $[0,2500000]$ (making sure that the highest coefficient at degree 8 is not 0) as a comparative basis: octic, since we will be contrasting with degree 8 Poincar\'e polynomials, upper bound of 2500000, since we can see from \fref{f:CY4-plot} and \eref{bCY4}, that the Min and Max of the Hodge numbers are respectively $(h^{1,1},h^{2,1},h^{3,1}) \in ([1,303148],[0,2010],[1,3030148])$, so that the upper bound to the $b_4$ term is 2425228.
Next, in part (b), we plot the same, but for monic palindromic octics, i.e., 
$P(t) = 1 + b_1 t + b_2 t^2 + b_3 t^3 + b_4 t^4 + b_3 t^5 + b_2 t^6 + b_1 t^7 + t^8$. 
Finally, in (c), we restrict once more, so that the linear term vanishes, i.e., $P(t) = 1 + b_2 t^2 + b_3 t^3 + b_4 t^4 + b_3 t^5 + b_2 t^6 + t^8$.

\begin{figure}[!h!t!b]
\centerline{
(a)
\includegraphics[trim=0mm 0mm 0mm 0mm, clip, width=3.5in]{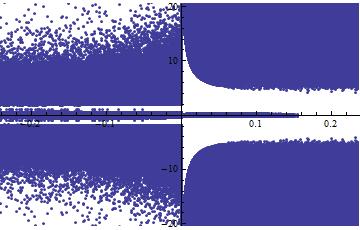}
(b)
\includegraphics[trim=0mm 0mm 0mm 0mm, clip, width=3.0in]{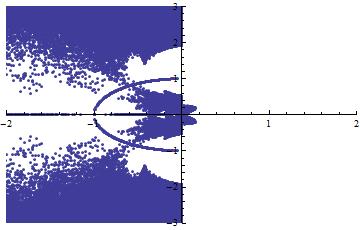}
}
\caption{{\sf 
(a) The position, on the complex plane, of the some 23 million zeros of the Poincar\'e polynomials of the approximately 1 million smooth Calabi-Yau fourfolds arising as hypersurfaces in toric five-folds.
(b) A slightly magnified area emphasizing the ordinate in the range $[-1,1]$.
}
\label{f:poincareCY4}}
\end{figure}

In contradistinction to these generic results, we can now find all the complex roots of the Poincar\'e polynomials of all known Calabi-Yau fourfolds.
The three million or so distinct Hodge data now presents a heavy computational challenge, on which a Quadra-core MacPro with 40Gb of memory laboured for a week, to produce some 23 million complex roots.
We present a scatter plot of these roots in part (a) of \fref{f:poincareCY4}.
In part (b) of the same figure, we magnify it slightly to emphasize the same range as the random plots in \fref{f:palin4}.

\begin{figure}[!h!t!b]
\centerline{
(a) \includegraphics[trim=0mm 0mm 0mm 0mm, clip, width=1.7in]{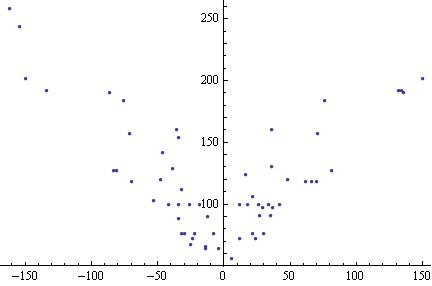}
(b) \includegraphics[trim=0mm 0mm 0mm 0mm, clip, width=1.7in]{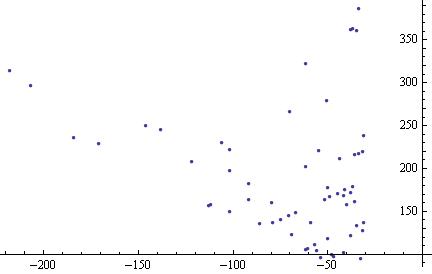}
(c) \includegraphics[trim=0mm 0mm 0mm 0mm, clip, width=2.5in]{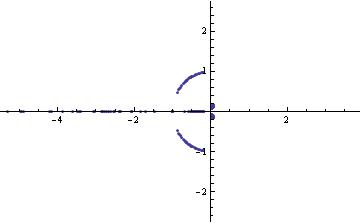}
}
\caption{{\sf 
(a) The Hodge numbers, with $h^{1,1}-h^{3,1}$ as the abscissa and $h^{1,1}+h^{3,1}$, the ordinate, of the fourfolds which have vanishing Euler number, and hence rank exceeding unity.
(b) The same, but with $h^{1,1}-h^{2,1}$ as the abscissa and $h^{1,1}+h^{2,1}$ as the ordinate.
(c) The position of the roots, on the complex plane, of the Poincar\'e polynomial of these 61 spaces out of the some 3 million.
}
\label{f:CY4chi=0}}
\end{figure}

\begin{figure}[!h!t!b]
\centerline{
(a) \includegraphics[trim=0mm 0mm 0mm 0mm, clip, width=2.5in]{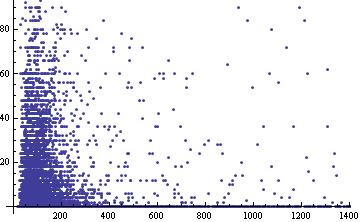}
(b) \includegraphics[trim=0mm 0mm 0mm 0mm, clip, width=3.5in]{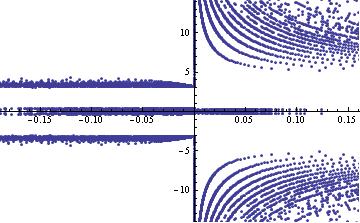}
}
\caption{{\sf 
(a) The Hodge numbers, with $h^{1,1}=h^{3,1}$ as the abscissa and $h^{2,1}$, the ordinate, of the ``self-mirror'' Calabi-Yau fourfolds.
(b) The position of the roots of the Poincar\'e polynomials of these 5009 spaces.
}
\label{f:selfmirrorCY4}}
\end{figure}

The theorem that $-1$ being a root of the Poincar\'e polynomial of $X$ implies $X$ has rank exceeding unity is generally applicable.
Hence, we can continue with this analysis.
Now, \eref{PCY4} implies that\footnote{Of course, the last equality follows directly for the definition of the the Euler characteristic $\chi(X)$.}
\begin{equation}\label{chi=0CY4}
P(-1,X) = 48 + 6 (h^{1,1} - h^{2,1} + h^{3,1})
= \chi(X) \ .
\end{equation}
However, the relation between $\chi$ and being self-mirror, or even the concept of the latter, is obviously not as clear in complex dimension greater than three.
Be that as it may, we can still examine \eref{chi=0CY4} in the the space of fourfolds.
Of the some 3 million distinct triplets, there are only 61 with vanishing Euler number, which we demonstrate in \fref{f:CY4chi=0}: in part (a), $h^{1,1}+h^{3,1}$ against $h^{1,1}-h^{3,1}$, and in part (b), $h^{1,1}+h^{2,1}$ against $h^{1,1}-h^{2,1}$.
We see that these are all of relatively small Hodge numbers, and in part (b), we see that in spite of the general linear behaviour seen in part (b) of \fref{f:CY4-plot}, there is some sub-structure.
In part (c), we plot the interesting shape of the roots of the Poincar\'e polynomials for these 61 members.
In retrospect to \fref{f:CY4-plot}, we see that perhaps the closest notion to mirror symmetry in complex dimension four is the interchange of $h^{1,1}$ and $h^{3,1}$. 
Of the circa 3 million, we find 5009 which have the property that $h^{1,1}=h^{3,1}$.
We plot the pairs $(h^{1,1}=h^{3,1}, h^{2,1})$ in part (a) of \fref{f:selfmirrorCY4} and the roots of their Poincar\'e polynomials in part (b).

\section{Non-Compact Calabi-Yau Geometries, Toric Diagrams and Newton Polynomials}
\setall
Having indulged in an excursion into the space of compact smooth Calabi-Yau threefolds and fourfolds, as well as their Poincar\'e polynomials, proceeding to the space of non-compact Calabi-Yau geometries is almost a perfunctory next step.
These are affine varieties such as flat space $\IC^d$ and singularities which locally admit Gorenstein resolutions, and are central to McKay Correspondence and generalizations in mathematics as well as AdS/CFT and branes in string theory.
A rich tapestry on this subject has been woven over the past few decades, whereby augmenting the relevance of our present investigation.

The most important family of non-compact Calabi-Yau geometries is indubitably those which afford toric description, as mentioned in the introduction.
In complex dimension three, the Calabi-Yau condition compels the toric diagram to be co-planar, whence each is characterized by a (convex) lattice polygon
$D = \{ v_i \}$, with each $v_i \in \IZ^2$.
Therefore, a polynomial which instantly springs to mind is the Newton polynomial\begin{equation}\label{Newton}
D = \{(x_i, y_i) \} \qquad
\Rightarrow \qquad 
P(z,w; X) = \sum\limits_i a_i z^{x_i} w^{y_i} \in \IC[z,w] \ ,
\end{equation}
where we have inserted potential coefficients $a_i$ for generality.
This is not a frivolous act;
indeed, when $a_i \in \IZ_{\ge 0}$, they are the so-called ``multiplicities'' first defined in \cite{Feng:2000mi} and play a vital r\^{o}le in comprehending the dimer model/brane tiling interpretation of toric gauge theories \cite{Hanany:2005ve,Franco:2005sm,Feng:2005gw}.

\begin{figure}[!h!t!b]
\centerline{\includegraphics[trim=0mm 0mm 0mm 0mm, clip, width=6.0in]{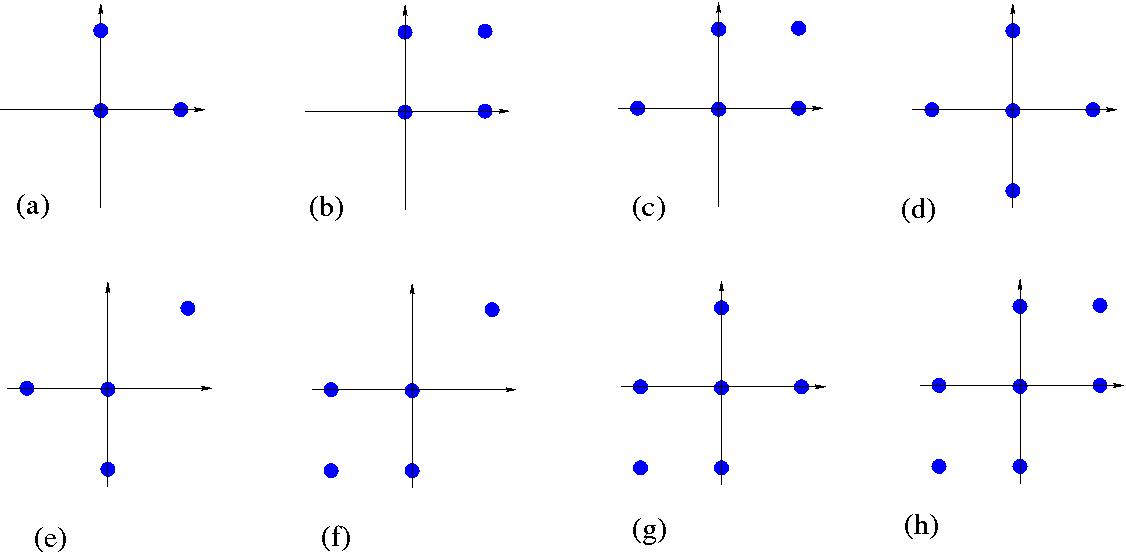}}
\caption{{\sf The most popular affine Calabi-Yau threefold toric diagrams, corresponding respectively to (a) $\IC^3$, (b) conifold, (c) SPP, (d) $\IF_0$, and (e,f,g,h) $dP_{0,1,2,3}$. The end points are at the standard lattice points in $\IZ^2$.}
\label{f:toricdiag}}
\end{figure}

The most famous toric diagrams for affine Calabi-Yau threefolds are depicted in \fref{f:toricdiag}, with the endpoints at the self-explanatory lattice points in $\IZ^2$; these include the reflexive polytopes in dimension two, and are commonly known as (a) $\IC^3$, (b) the conifold, (c) the suspended pinched point (SPP), (d) affine cone over the zeroth Hirzebruch surface $\IF_0 = \IP^1 \times \IP^1$, (e) $dP_0 = \IC^3 / \IZ_3$, the affine cone over $\IP^2$, and (f,g,h) $dP_n$, cones over respectively the first, second and third del Pezzo surfaces which are $\IP^2$ blown up at $n=$ 1,2,and 3 generic points.

\begin{figure}[!h!t!b]
\centerline{
\includegraphics[trim=0mm 0mm 0mm 0mm, clip, width=1.7in]{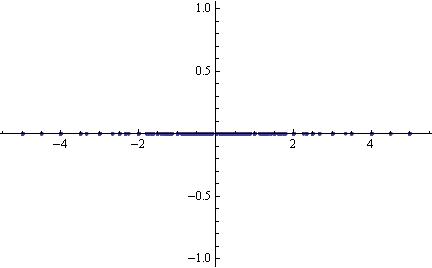}
\includegraphics[trim=0mm 0mm 0mm 0mm, clip, width=1.7in]{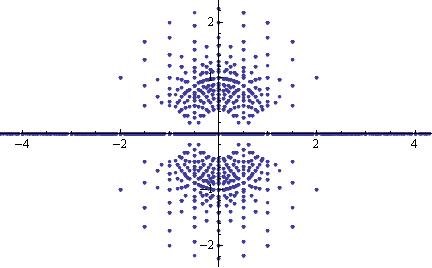}
\includegraphics[trim=0mm 0mm 0mm 0mm, clip, width=1.7in]{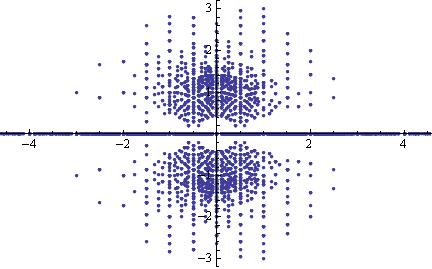}
\includegraphics[trim=0mm 0mm 0mm 0mm, clip, width=1.7in]{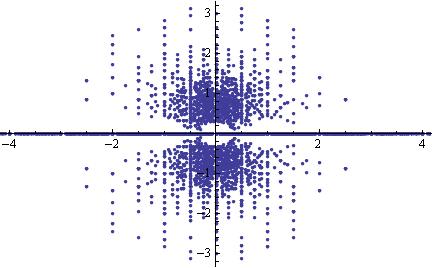}
}
\caption{{\sf The roots of the Newton polynomials $P(z,w)$ at $z=1$, for (a) the conifold, (b) $\IF_0$, (c) $dP_1$ and (d) $dP_3$. 
}
\label{f:proj1}}
\end{figure}

An immediate difficulty with \eref{Newton} is, of course, that the polynomial is bi-variate, whereby describing, algebraically, Riemann surfaces. Even though such surfaces are crucial in the understanding of the gauge theory constructed on branes probing these affine toric Calabi-Yau spaces\footnote{Cf.~\cite{Feng:2005gw} for discussion on a web of inter-relations and various projections and spines of these Riemann surfaces.}, the notion of zeros is not obvious.
We could, for example, set one of the co\"ordinates to a fixed value, and consider the roots of the resulting uni-variate projection.
This, however, does not seem particularly natural.
Nevertheless, for illustrative purposes, we include a few examples in \fref{f:proj1}, wherein we have set $z$ to 1, varied the coefficients $a_i$ randomly and integrally in $[-5,5]$, and plotted the roots of the resulting polynomial in $w$ for 5000 samples.

A much more natural and, as it turns out, interesting direction to take is to consider \eref{Newton} not as a complex, but as a real, curve.
For comparative purposes, we should be mindful of the ``amoebae'' and so-dubbed ``algae'' projections proposed in \cite{Feng:2005gw}, which are, respectively the (natural log of the) real and imaginary projections of the Newton polynomials of the associated toric Calabi-Yau threefold.

In this regard, perhaps the most significant quantity for our complex Newton polynomials is the set of turning points in $\IR^2$, viz., the set of real critical points of $P(z,w) \in \IR[z,w]$.
In other words, we find the simultaneous real solutions to
\begin{equation}\label{critpt}
\partial_z P(z,w) = \sum\limits_i a_i x_i z^{x_i-1} w^{y_i} = 0 \ , \quad
\partial_w P(z,w) = \sum\limits_i a_i y_i z^{x_i} w^{y_i-1} = 0 \ ,
\end{equation}
with $a_i$ randomly sampled and discard the imaginary solutions.
This, indeed, brings us back to Waring's original considerations on non-real roots.

\begin{figure}[!h!t!b]
\centerline{
\begin{tabular}{c}
(b) \includegraphics[trim=0mm 0mm 0mm 0mm, clip, width=1.7in]{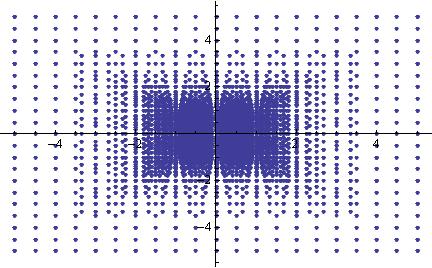}
(c) \includegraphics[trim=0mm 0mm 0mm 0mm, clip, width=1.7in]{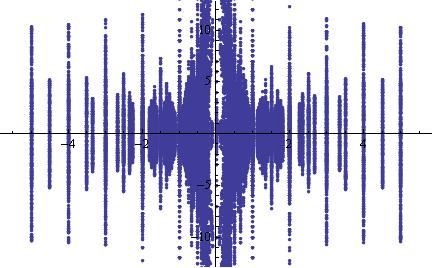}
(d) \includegraphics[trim=0mm 0mm 0mm 0mm, clip, width=1.7in]{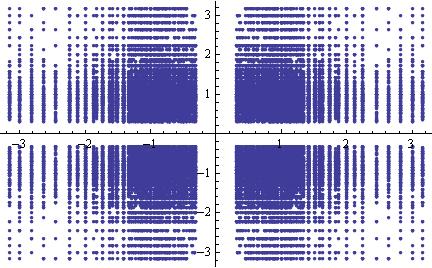}
\\
~\\
(e) \includegraphics[trim=0mm 0mm 0mm 0mm, clip, width=1.7in]{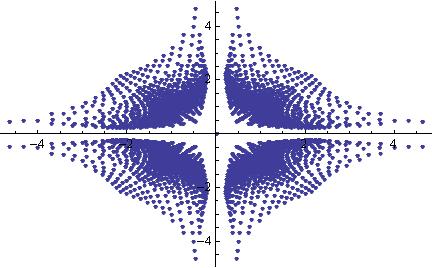}
(f) \includegraphics[trim=0mm 0mm 0mm 0mm, clip, width=1.7in]{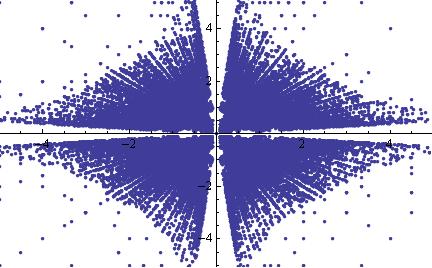}
(g) \includegraphics[trim=0mm 0mm 0mm 0mm, clip, width=1.7in]{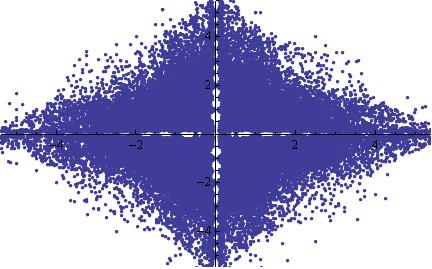}
(h) \includegraphics[trim=0mm 0mm 0mm 0mm, clip, width=1.7in]{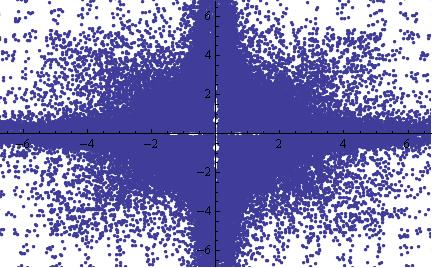}
\end{tabular}
}
\caption{{\sf In reference to the toric diagrams in \fref{f:toricdiag}, for each Calabi-Yau threefold geometry, we find the corresponding Newton polynomial as a real algebraic curve. We sample over 50000 random integer coefficients in the range $[-10,10]$, and isolate the real critical points in $\IR^2$, which are then plotted collectively for (b) conifold, (c) SPP, (d) $\IF_0$, and (e,f,g,h) $dP_{0,1,2,3}$.
}
\label{f:critical}}
\end{figure}

In \fref{f:critical}, we take a fixed toric diagram corresponding to a given Calabi-Yau geometry, and consider the Newton Polynomial in \eref{Newton}.
Then, we sample 50000 random integer coefficients $a_i$ in an appropriate range, here taken to be $[-10,10]$. For each, we find the real critical points, and collectively plot them.
We note that $\IC^3$ does not that any real critical points and is thus left out.
This is because the Newton polynomial is simply $a+bz+cw$ for $a,b,c \in \IZ$.
Thus, \eref{critpt} gives $b=c=0$, independent of $(z,w)$ co\"odinates.
Similarly, the case of (b), the conifold, can be considered a reference point.
The Newton polynomial is $a + bz +cw + d z w$, whence, the critical points are given by the solutions of $b + d w = c + d z = 0$, or, $w_0 = -b/d, z_0 = -c/d$.
Hence, given that each of $b,c,d$ is independently randomly evenly distributed, the turning points $(z_0,w_0)$ are then distributed as quotients of even random samples, and whence the clustering nearer to the lower values as seen in the darker region in the figure.

\section*{Acknowledgements}
{\it 
Ad Catharinae Sanctae Alexandriae et ad Maiorem Dei Gloriam, cum Universitate Nankai, Universitate Civitate Londiniensis, Scientiae et Technologiae Concilio Anglicae, et Collegio Mertonense Oxoniensis, atque amore Elizabetae Katherinae Hunter, Rosae Hiberniae, multas gratias YHH agit.
}

{\small
\begin{singlespace}

\end{singlespace}
}
\end{document}